\documentclass[11pt,twoside]{article}
\usepackage[utf8]{inputenc}
\usepackage[english]{babel}
\usepackage{fancyhdr}
\usepackage{amssymb,amsmath,amsthm}
\usepackage{blindtext}
\usepackage{hyperref}
\usepackage{caption}
\usepackage[square,numbers]{natbib}
\usepackage{graphicx}
\usepackage{amsmath}
\usepackage{amssymb}
\usepackage{latexsym}
\usepackage{algorithmicx}
\usepackage{algorithm}
\usepackage{algpseudocode}
\usepackage{subfig}
\usepackage{epsfig}
\usepackage{longtable}
\usepackage[inline]{enumitem}
\usepackage{xcolor}

\newtheorem{proposition}{Proposition}
\newtheorem{definition}{Definition}
\newtheorem{lemma}{Lemma}
\newtheorem{theorem}{Theorem}

\def\R{\mathbb R}

\def\eps{\epsilon}
\def\beq{\begin{equation}}
\def\eeq{\end{equation}}

\def\Re{\mathbb{R}}
\def\x-{{\bar x}}

\def\dg+{(\nabla g_i(\x-)'d)^{\scriptscriptstyle +}}

\def\to{\rightarrow}

\hypersetup{
    colorlinks=false,
    urlcolor=black,
    pdfborder={0 0 0}
}

\newcommand{\email}[1]{E-mail: \href{mailto:#1}{\texttt{#1}}}

\begin{document}
\thispagestyle{plain}

\setcounter{page}{1}

{\centering
{\LARGE \bfseries An augmented Lagrangian method exploiting an active-set strategy and second-order information}

\bigskip\bigskip
Andrea Cristofari$^*$, Gianni Di Pillo$^\dag$, Giampaolo Liuzzi$^\dag$, Stefano Lucidi$^\dag$
\bigskip

}

\begin{center}
\small{\noindent$^*$Department of Mathematics ``Tullio Levi-Civita'' \\
University of Padua \\
Via Trieste, 63, 35121 Padua, Italy \\
\email{andrea.cristofari@unipd.it} \\
\bigskip
$^\dag$Department of Computer, Control and Management Engineering \\
Sapienza University of Rome \\
Via Ariosto, 25, 00185 Rome, Italy \\
E-mail: \texttt{gianni.dipillo@diag.uniroma1.it}, \texttt{giampaolo.liuzzi@diag.uniroma1.it}, \texttt{stefano.lucidi@diag.uniroma1.it} \\
}
\end{center}

\bigskip\par\bigskip\par
\noindent \textbf{Abstract.}
In this paper, we consider nonlinear optimization problems with nonlinear equality constraints and bound constraints on the variables. For the solution of such problems, many augmented Lagrangian methods have been defined in the literature. Here, we propose to modify one of these algorithms, namely ALGENCAN by Andreani et al., in such a way to incorporate second-order information into the augmented Lagrangian framework, using an active-set strategy. We show that the overall algorithm has the same convergence properties as ALGENCAN and an asymptotic quadratic convergence rate under suitable assumptions. The numerical results confirm that the proposed algorithm is a viable alternative to ALGENCAN with greater robustness.

\bigskip\par
\noindent \textbf{Keywords.} Constrained optimization. Augmented Lagrangian methods. Nonlinear programming algorithms. Large-scale optimization.

\bigskip\par
\noindent \textbf{MSC2000 subject classifications.} 90C30. 65K05.

\section{Introduction}\label{sec:intro}

In this paper, we are interested in the solution of smooth constrained
optimization problems  of the type:
\begin{equation}\label{probP}
\begin{array}{lll}
&\min & f(x)\\&& h(x) = 0 \\ && \ell \leq x \leq u,
\end{array}\end{equation}
where $x,\ell,u\in \Re^n,$ $\ell_i < u_i$, for all $i=1,\dots,n$, $f \colon \Re^n\to \Re$, $h \colon \Re^n\to \Re^p$
are twice continuously differentiable functions. Note that the structure of Problem~\eqref{probP} is sufficiently general to capture, through reformulation, also problems with nonlinear inequality constraints. Problem (\ref{probP}) has been studied for decades and many optimization methods have been proposed for its solution. Solution algorithms for~\eqref{probP} belong to different classes like, e.g., sequential penalty \cite{fiacco1990nonlinear}, augmented Lagrangian \cite{bertsekas2014constrained} and sequential quadratic programming \cite{nocedal2006numerical}.

Among the algorithms based on augmented Lagrangian functions, the one implemented in the ALGENCAN \cite{algencan:2008,andreani:2008} software package is one of the latest and more efficient. The computational heavy part of ALGENCAN consists in the solution (at every outer iteration) of the subproblem, i.e., the minimization of the augmented Lagrangian merit function for given values of the penalty parameter and of the estimated Lagrange multipliers. Such minimization is carried out by the inner solver GENCAN~\cite{birgin2002large}.

It is worth noticing that besides the above methods, efficient local algorithms have been proposed in the literature that exploit second-order information to define superlinearly convergent Newton-like methods~\cite{bertsekas2014constrained,dipillo:2000,facchinei:1995}. 
The so-called ``acceleration strategy'' of ALGENCAN is an attempt to exploit second-order information by means of such locally convergent methods to improve the convergence rate of the overall algorithm.

The idea that we develop in this paper is twofold. On the one side, we propose an alternative and possibly more extensive way to use second-order information within the framework of an augmented Lagrangian algorithm. Basically, we propose a Newton-type direction to use even when potentially far away from solution points.
The use of such a Newton direction is combined with an appropriate active-set strategy.
In particular, after estimating active and non-active variables with respect to the bound constraints,
we compute the Newton direction with respect to only the variables estimated as non-active,
while the ones estimated as active are set to the bounds.

On the other hand, when the Newton-type direction cannot be computed or does not satisfy a proper condition, we propose to resort to the minimization of the augmented Lagrangian function, but using an efficient active-set method for bound-constrained problems~\cite{cristofari:2017}.

The paper is organized as follows. In Section~\ref{sec:preliminary}, we report some preliminary results that will be useful in the paper. In Section~\ref{sec:direction}, we describe the procedure to compute the Newton-type direction and we study its theoretical properties. Section~\ref{sec:algorithm} is devoted to the description of the proposed augmented Lagrangian algorithm and to its convergence analysis. In Section~\ref{sec:rate}, we are concerned with the analysis of the converge rate for the proposed method. In Section~\ref{sec:experiments}, we report some numerical experiments and comparison with existing software. Finally, in Section~\ref{sec:conclusions} we draw some conclusions.

\section{Notation and Preliminary Results}\label{sec:preliminary}
Given a vector $x \in \R^n$, we denote by $x_i$ its $i$th entry and, given an index set $T \subseteq \{1,\ldots,n\}$, we denote by $x_T$ the subvector obtained from $x$ by discarding the components not belonging to $T$.
The gradient of a function $f(x)$ is denoted by $\nabla f(x)$, while the Hessian matrix is denoted by $\nabla^2 f(x)$.
We indicate by $\nabla_{x_i} f(x)$ the $i$th entry of $\nabla f(x)$.
The Euclidean norm of a vector $x$ is indicated by $\|x\|$, while $\|x\|_{\infty}$ denotes the sup-norm of $x$.
Given a matrix $M$, we indicate by $\|M\|$ the matrix norm induced by the Euclidean vector norm.
The projection of a vector $x$ onto a box $[a,b]$ is denoted by ${\cal P}_{[a,b]}(x)$.
The $i$th column of the identity matrix is indicated by $e_i$.

With reference to Problem~\eqref{probP},
we define the Lagrangian function $L(x,\mu)$ with respect to the equality constraints as follows:
\[
L(x,\mu) := f(x)+\mu^T h(x),
\]
where $\mu\in\Re^p$ is the Lagrange multiplier.

Denoting the gradient of $L(x,\mu)$ with respect to $x$ as $\nabla_x L(x,\mu) = \nabla f(x) + \nabla h(x) \mu$,
we say that $(x,\mu,\sigma,\rho)\in \Re^{3n+p}$ is a KKT tuple for Problem~\eqref{probP} if
\begin{subequations}\label{KKT_conditions}
\begin{align}
\nabla_x L(x,\mu) & = \sigma - \rho, \label{KKT_conditions_1}\\
h(x) & = 0, \\
\sigma^T (\ell - x) & = 0, \label{KKT_conditions_3} \\
\rho^T (x-u) & = 0, \label{KKT_conditions_4} \\
\ell - x & \le 0,\quad \sigma \ge 0, \\
x - u & \le 0, \quad \rho \ge 0.
\end{align}
\end{subequations}
If $x^*$ is local minimum of Problem~\eqref{probP} that satisfies some constraint qualification, then there exist KKT multipliers $\mu^*,\sigma^*,\rho^*$ such that
$(x^*,\mu^*,\sigma^*,\rho^*)$ is a KKT tuple.
Note that the KKT conditions~\eqref{KKT_conditions} can be rewritten as follows:
\begin{subequations}\label{KKT_conditions2}
\begin{align}
\nabla_{x_i} L(x,\mu) &
\begin{cases}
= 0, \quad & \text{if } \ell_i < x_i < u_i, \\
\ge 0, \quad & \text{if } x_i = \ell_i, \\
\le 0, \quad & \text{if } x_i = u_i,
\end{cases} \\
h(x) & = 0.
\end{align}
\end{subequations}

For a KKT tuple $(x^*,\mu^*,\sigma^*,\rho^*)$, we say that the \textit{strict complementarity} holds if
$x^*_i = \ell_i \Rightarrow \sigma^*_i > 0$ and $x^*_i = u_i \Rightarrow \rho^*_i > 0$, that is,
$x^*_i = \ell_i \Rightarrow \nabla_i L(x^*,\mu^*) > 0$ and $x^*_i = u_i \Rightarrow \nabla_i L(x^*,\mu^*) < 0$.


Now, let us define the multiplier functions $\sigma(x,\mu)$ and $\rho(x,\mu)$,
which give us some estimates of the KKT multipliers $\sigma$ and $\rho$, respectively, associated to the box constraints of Problem~\eqref{probP}.
Following the same approach used in~\cite{cristofari:2017,desantis2012active} for bound-constrained problems, we can first express $\sigma(x,\mu) = \nabla_x L(x,\mu) + \rho(x,\mu)$ from~\eqref{KKT_conditions_1}, and then we can compute $\rho(x,\mu)$ by minimizing the error over~\eqref{KKT_conditions_3}--\eqref{KKT_conditions_4}
(see~\cite{desantis2012active} for more details), obtaining
\begin{align}
\label{def:sigma}\sigma_i(x,\mu)  & :=  \phantom{-}\frac{(u_i-x_i)^2}{(\ell_i-x_i)^2 + (u_i-x_i)^2} \nabla_{x_i}L(x,\mu), \quad i = 1,\ldots,n, \\
\label{def:rho}\rho_i(x,\mu) & :=  -\frac{(\ell_i-x_i)^2}{(\ell_i-x_i)^2 + (u_i-x_i)^2} \nabla_{x_i}L(x,\mu), \quad i = 1,\ldots,n.
\end{align}
These multiplier functions will be employed later for defining an active-set strategy to be used in the proposed algorithm.

Moreover, now we can say that
$(x^*,\mu^*)\in\Re^{n+p}$ is a KKT pair for Problem~\eqref{probP} when $(x^*,\mu^*,\sigma(x^*,\mu^*),\rho(x^*,\mu^*))$ is a KKT tuple.

\subsection{The Augmented Lagrangian Method}\label{subsec:auglag}
The algorithm we propose here builds upon the augmented Lagrangian method described in~\cite{andreani:2008}, where an augmented Lagrangian function is defined with respect to a subset of constraints and iteratively minimized over $x$ subject to the remaining constraints.
In our case, we define the augmented Lagrangian function for Problem~\eqref{probP} with respect to the equality constraints as
\[
L_a(x,\mu;\eps) := L(x,\mu) + \frac{1}{\eps}\|h(x)\|^2,
\]
where $\eps > 0$ is a parameter that penalizes violation of the equality constraints.
Given an estimate $(x_k,\bar \mu_k)$ of a KKT pair and a value $\epsilon_k$ for the penalty parameter, the new iterate $x_{k+1}$ can thus be computed by approximately solving the following bound-constrained subproblem:
\begin{equation}\label{subprobP}
\begin{array}{l}
\min \ L_a(x,\bar\mu_k;\eps_k)\\~~~~  \ell \leq x \leq u.
\end{array}\end{equation}
Then, according to~\cite{andreani:2008}, we can set
\begin{equation}\label{update_molt}
 \mu_{k+1} = \bar \mu_k + \frac{2}{\eps_k}h(x_{k+1})
\end{equation}
and update the Lagrange multiplier $\bar \mu_{k+1}$ by projecting $(\mu_{k+1})_i$ in a suitable interval $[\bar\mu_{\text{min}},\bar\mu_{\text{max}}]$, $i = 1,\ldots,p$,
that is,
\begin{equation}\label{update_molt_2}
(\bar \mu_{k+1})_i = \max\{\bar\mu_{\text{min}},\min\{(\mu_{k+1})_i,\bar\mu_{\text{max}}\}\}, \quad i=1,\ldots,p.
\end{equation}
Finally, we decrease the penalty parameter $\epsilon_{k+1}$ if the constraint violation is not sufficiently reduced and start a new iteration.
We can summarize the method proposed in~\cite{andreani:2008} as in the following scheme.

\begin{algorithm}[h]
\caption*{\centering Augmented Lagrangian Method}
\begin{algorithmic}
\item[] {\bf Given} finite scalars $\bar\mu_{\text{min}} < \bar\mu_{\text{max}}$, $\beta \in (0,1)$, $\eta \in (0,1)$, $\theta \in (0,1)$, $\eps_0>0$, a sequence \mbox{$\{\tau_k\}\searrow 0$},  a starting point $x_0\in [\ell,u]$ and estimates of multipliers $(\bar \mu_0)_i = (\mu_0)_i \in [\bar\mu_{\text{min}},\bar\mu_{\text{max}}]$, $i=1,\ldots,p$

 \item[] {\bf For} $k=0,1,\dots$
 \item[]\makebox[0.7cm][l]{}{\bf Compute} $x_{k+1}$ as an approximate solution of~\eqref{subprobP} with tolerance $\tau_k$

      \item[]\makebox[0.7cm][l]{}{\bf Set} $\mu_{k+1}$ by~\eqref{update_molt} and $\bar \mu_{k+1}$ by~\eqref{update_molt_2}

      \item[]\makebox[0.7cm][l]{}{\bf If} $\|h(x_{k+1})\|_{\infty} \le \eta \|h(x_k)\|_{\infty}$, then set $\eps_{k+1} = \eps_k$, else set $\eps_{k+1} = \theta \eps_k$
 \item[] {\bf End for}
\end{algorithmic}
\end{algorithm}

In the next section, we will describe how to incorporate the use of a proper second-order direction into this augmented Lagrangian framework.

\section{Direction Computation}\label{sec:direction}
In this section we introduce and analyze the procedure for computing a second-order direction, employing a proper active-set estimate.

\subsection{Active-Set Estimate}
Taking inspiration from the strategy proposed in~\cite{facchinei:1995},
for any $x \in [\ell,u]$ and any $\mu \in \R^p$,
we can estimate the active constraints in a KKT point by the following sets:
\begin{gather}
{\cal L}(x,\mu) :=  \{i \colon \nabla_{x_i} L(x,\mu) > 0,\ \ell_i \leq x_i \leq \ell_i+ \nu\sigma_i(x,\mu)\}, \label{L_as}\\
{\cal U}(x,\mu) := \{i \colon \nabla_{x_i} L(x,\mu) < 0,\ u_i - \nu\rho_i(x,\mu) \leq x_i \leq u_i\}, \label{U_as}
\end{gather}
where $\nu>0$ is a given parameter and
the multiplier functions $\sigma(x,\mu)$, $\rho(x,\mu)$ are defined in~\eqref{def:sigma} and~\eqref{def:rho}, respectively.

In particular, in a given pair $(x,\mu)$, the sets ${\cal L}(x,\mu)$ and ${\cal U}(x,\mu)$ contain the indices of the variables that are estimated to be active at the lower bound $\ell_i$ and at the upper bound $u_i$, respectively, in a KKT point.
As to be shown later, at each iteration of the proposed algorithm, these sets are used to compute a Newton direction with respect to only the variables that are estimated as non-active, while the variables estimated as active are set to bound.

Using results from~\cite{facchinei:1995}, the following identification property of the active-set estimate~\eqref{L_as}--\eqref{U_as} holds.
\begin{proposition}\label{prop:estim}
If $(x^*,\mu^*,\sigma^*,\rho^*)$ satisfies the KKT conditions~\eqref{KKT_conditions}, then there exists a neighborhood of $(x^*,\mu^*)$ such that, for each $(x,\mu)$ in this neighborhood, we have
\begin{gather*}
\{i \colon x^*_i = \ell_i, \, \sigma^*_i > 0 \} \subseteq {\cal L}(x,\mu) \subseteq \{i \colon x^*_i = l_i\}, \\
\{i \colon x^*_i = u_i, \, \rho^*_i > 0 \} \subseteq {\cal U}(x,\mu) \subseteq \{i \colon x^*_i = u_i\}.
\end{gather*}

In particular, if the strict complementarity holds at $(x^*,\mu^*,\sigma^*,\rho^*)$, for each $(x,\mu)$ in this neighborhood we have
\[
{\cal L}(x,\mu) = \{i \colon x^*_i = l_i\} \quad \text{and} \quad
{\cal U}(x,\mu)  = \{i \colon x^*_i = u_i\}.
\]
\end{proposition}

The result stated in the above proposition holds for an unknown neighborhood of the optimal solution. It would be of great interest and importance to give a characterization of that neighborhood, in order to bound the maximum number of iterations required by the algorithm to identify the active set.
Currently, this is an open problem and we think it may represent a possible line of future research, for example by adapting the complexity results given for ALGENCAN in~\cite{birgin2020complexity}, or extending some results on finite active-set identification given in the literature for specific classes of algorithms~\cite{bomze:2019b,cristofari2021active,nutini:2019}.

\subsection{Step Computation}
In the proposed algorithm, at the beginning of every iteration $k$, we have a point $x_k \in [\ell,u]$ and Lagrange multiplier estimates $(\bar \mu_k)_i \in [\bar\mu_{\text{min}},\bar\mu_{\text{max}}]$, $i=1,\ldots,p$.

Using~\eqref{L_as}--\eqref{U_as}, we estimate the active and non-active set in $(x_k,\bar \mu_k)$. Denoting
\begin{equation}\label{active_set_estimate_k}
{\cal L}_k := {\cal L}(x_k,\bar \mu_k), \quad {\cal U}_k := {\cal U}(x_k, \bar\mu_k), \quad {\cal B}_k := {\cal L}_k\cup{\cal U}_k,
\quad {\cal N}_k := \{1,\ldots,n\} \setminus {\cal B}_k,
\end{equation}
we can thus partition the vector $x_k$ as $x_k = (x_{{\cal B}_k}, x_{{\cal N}_k})$, reordering its entries if necessary.
Let us also denote
\[
L_k := L(x_k,\bar \mu_k), \quad \nabla_{{\cal N}_k} L_k := [\nabla_x L_k]_{{\cal N}_k}, \quad h_k := h(x_k), \quad \nabla_{{\cal N}_k} h_k := [\nabla h_k]_{{\cal N}_k},
\]
while $\nabla^2_{xx} L_k$ denotes the Hessian matrix of $L_k$ deriving with respect to $x$ two times and $\nabla^2_{{\cal N}_k} L_k$ denotes the submatrix obtained from $\nabla^2_{xx} L_k$ by discarding rows and columns not belonging to ${\cal N}_k$.

Now, consider the following system of equation with unknowns $x_{{\cal N}_k}$ and $\mu$:
\begin{subequations}\label{KKT_red}
\begin{align}
\nabla_{{\cal N}_k} L(x_{{\cal N}_k}, x_{{\cal B}_k},\mu) & = 0,\label{KKT_reda}\\
h(x_{{\cal N}_k}, x_{{\cal B}_k}) & = 0.\label{KKT_redb}
\end{align}
\end{subequations}
The nonlinear system \eqref{KKT_reda}--\eqref{KKT_redb} can be solved iteratively by the Newton method, where the Newton direction is computed by
solving the following linear system:
\begin{equation}\label{KKT_systemB}
\begin{pmatrix}\nabla^2_{{\cal N}_k} L_k & \nabla_{{\cal N}_k} h_k \\ \nabla_{{\cal N}_k} h_k^T & 0\end{pmatrix}
 \begin{pmatrix}d_{x_{{\cal N}_k}} \\ d_{\mu}\end{pmatrix}
 = -\begin{pmatrix}\nabla_{{\cal N}_k}L_k \\ h_k\end{pmatrix}.
\end{equation}
Hence, if a solution $(d_{x_{{\cal N}_k}},d_{\mu})$ of~\eqref{KKT_systemB} exists, we can set
\[
d_k = (d_{x_{{\cal N}_k}},d_{\mu})
\]
and move from $((x_k)_{{\cal N}_k},\bar \mu_k)$ along $d_k$,
then projecting $(x_k)_{{\cal N}_k} + d_{x_{{\cal N}_k}}$ onto the box $[\ell_{{\cal N}_k},u_{{\cal N}_k}]$.
In particular, we define
\begin{equation}\label{set_non_active_vars}
(\tilde x_k)_{{\cal N}_k} = {\cal P}_{[\ell_{{\cal N}_k},u_{{\cal N}_k}]}((x_k)_{{\cal N}_k} + d_{x_{{\cal N}_k}}).
\end{equation}
and
\[
\mu_{k+1} = \bar \mu_k + d_{\mu}.
\]
For what concerns the variables $(x_k)_{{\cal B}_k}$,
since they are estimated as active, we set them to the bounds.
Namely, we define $(\tilde x_k)_{{\cal B}_k}$ as follows:
\begin{equation}\label{set_active_vars}
(\tilde x_k)_i =
\begin{cases}
\ell_i, \quad & \text{if } i\in {\cal L}_k, \\
u_i, \quad & \text{if } i\in {\cal U}_k.
\end{cases}
\end{equation}
The following results holds.

\begin{proposition}
If the  solution $d_k$ of system~\eqref{KKT_systemB} exists,
then $(x_k,\bar \mu_k,\sigma_k,\rho_k)$ is a KKT tuple
with $\sigma_k = \sigma(x_k,\bar \mu_k)$ and $\rho_k = \rho(x_k,\bar \mu_k)$ if and only if
$d_k=0$ and $(\tilde x_k)_{{\cal B}_k} = (x_k)_{{\cal B}_k}$.
\end{proposition}

{\it Proof} First, assume that $d_k=0$ and $(\tilde x_k)_{{\cal B}_k} = (x_k)_{{\cal B}_k}$. From~\eqref{KKT_systemB}, we have
\[
\nabla_{x_{{\cal N}_k}} L(x_k,\bar \mu_k) = 0 \quad \text{and} \quad h(x_k) = 0.
\]
Using the expression of ${\cal L}(x_k,\bar \mu_k)$ and ${\cal U}(x_k,\bar \mu_k)$ given in~\eqref{L_as}--\eqref{U_as}, and recalling the definition of $\rho(x,\mu)$ and $\sigma(x,\mu)$ given in~\eqref{def:sigma}--\eqref{def:rho}, we also have
\begin{align*}
(\sigma_k)_i = (\rho_k)_i = \nabla_{x_i} L(x_k,\bar \mu_k) = 0, \qquad & \forall i \in {\cal N}_k, \\
(x_k)_i = (\tilde x_k)_i = \ell_i, \quad (\sigma_k)_i = \nabla_{x_i} L(x_k,\bar \mu_k) > 0, \quad (\rho_k)_i = 0, \qquad & \forall i \in {\cal L}_k, \\
(x_k)_i = (\tilde x_k)_i = u_i, \quad (\sigma_k)_i = 0, \quad (\rho_k)_i = -\nabla_{x_i} L(x_k,\bar \mu_k) > 0, \qquad & \forall i \in {\cal U}_k.
\end{align*}
It follows that KKT conditions~\eqref{KKT_conditions} are satisfied.

Now, assume that $(x_k,\bar \mu_k,\sigma_k,\rho_k)$ is a KKT tuple.
Since $\nabla_{x_{{\cal N}_k}} L(x_k,\bar \mu_k) = 0$ and $h(x_k) = 0$, from~\eqref{KKT_systemB} we have $d_k=0$. Finally, using the KKT conditions written as in~\eqref{KKT_conditions2}, and recalling the definition of $\rho(x,\mu)$ and $\sigma(x,\mu)$ given in~\eqref{def:sigma}--\eqref{def:rho}, we also have $(x_k)_i = \ell_i = (\tilde x_k)_i$ for all $i \in {\cal L}_k$ and $(x_k)_i = u_i = (\tilde x_k)_i$ for all $i \in {\cal U}_k$.
\qed

\section{The Algorithm}\label{sec:algorithm}
In this section, we use the above described active-set estimate and Newton strategy to design a primal-dual augmented Lagrangian method.

At the beginning of each iteration $k$, we have a pair $(x_k,\bar \mu_k)$. We first
estimate the active set ${\cal L}_k \cup {\cal U}_k$ and the non-active set ${\cal N}_k$ as in~\eqref{active_set_estimate_k}.
If possible, we calculate a direction $d_k = (d_{x_{\cal N}},d_{\mu})$ by solving the Newton system~\eqref{KKT_systemB} and we compute $(\tilde x_k)$ as in~\eqref{set_non_active_vars}--\eqref{set_active_vars}.
This point is accepted and set as $x_{k+1}$ only if $\|(d_k,(\tilde x_k-x_k)_{{\cal B}_k})\| \le \Delta_k$, where $\Delta_k$ is iteratively decreased trough the iterations by a factor $\beta \in (0,1)$.

If this is not the case, we compute $x_{k+1}$ as an approximate minimizer of the bound-constrained subproblem~\eqref{subprobP}, such that
\begin{equation}\label{opt_subprob}
\|x_{k+1} - {\cal P}_{[\ell,u]}(x_{k+1} - \nabla_x L_a(x_{k+1},\bar \mu_k;\eps_k))\|_{\infty} \le \tau_k,
\end{equation}
with $\{\tau_k\} \to 0$.
Then, we update the multiplier estimate $\mu_{k+1}$ by~\eqref{update_molt} and decrease the penalty parameter $\epsilon_{k+1}$ if the constraint violation is not sufficiently reduced.

We finally terminate the iteration by setting $\bar \mu_{k+1}$ as the projection of $\mu_{k+1}$ on a prefixed box, according to~\eqref{update_molt_2}.

The proposed method, named Primal-Dual Augmented Lagrangian Method (\mbox{P-D ALM}), is reported in the following algorithmic scheme.
As specified later (see Section~\ref{sec:experiments}), in practical implementation of the algorithm we use a stricter test to accept the point $\tilde x_k$, also requiring a decrease of the feasibility violation in the new point $\tilde x_k$.
For the sake of generality, the theoretical analysis is carried out by considering only the condition $\|(d_k,(\tilde x_k-x_k)_{{\cal B}_k})\| \le \Delta_k$.

\begin{algorithm}[h]
\caption*{\centering Primal-Dual Augmented Lagrangian Method (\mbox{P-D ALM})}
\label{alg:pdalm}
\begin{algorithmic}
\item[] {\bf Given} finite scalars $\bar\mu_{\text{min}} < \bar\mu_{\text{max}}$, $\beta \in (0,1)$, $\eta \in (0,1)$, $\theta \in (0,1)$, $\Delta_0>0$, $\eps_0>0$, a sequence \mbox{$\{\tau_k\}\searrow 0$},  a starting point $x_0\in [\ell,u]$ and estimates of multipliers $(\bar \mu_0)_i = (\mu_0)_i \in [\bar\mu_{\text{min}},\bar\mu_{\text{max}}]$, $i=1,\ldots,p$

 \item[] {\bf For} $k=0,1,\dots$
    \item[]\makebox[0.7cm][l]{}{\bf Compute} the active and non-active set estimates ${\cal L}_k, {\cal U}_k, {\cal N}_k$ as in~\eqref{active_set_estimate_k}

    \item[]\makebox[0.7cm][l]{}{\bf Compute} $d_k = (d_{x_{\cal N}},d_{\mu})$ by solving~\eqref{KKT_systemB}, if possible, and set $(\tilde x_k)$ as in~\eqref{set_non_active_vars}--\eqref{set_active_vars}

    \item[]\makebox[0.7cm][l]{}{\bf If} $d_k$ has been computed and $\|(d_k,(\tilde x_k-x_k)_{{\cal B}_k})\| \le \Delta_k$, then set $x_{k+1} = \tilde x_k$,\\
    \makebox[1.1cm][l]{}$\mu_{k+1} = \bar \mu_k + d_{\mu}$, $\Delta_{k+1} = \beta \Delta_k$ and $\eps_{k+1} = \eps_k$

    \item[]\makebox[0.7cm][l]{}{\bf Else}, compute $x_{k+1}$ satisfying~\eqref{opt_subprob} and set $\mu_{k+1}$ by~\eqref{update_molt}.\\
    \makebox[1.1cm][l]{}{\bf If} $\|h(x_{k+1})\|_{\infty} \le \eta \|h(x_k)\|_{\infty}$, then set $\eps_{k+1} = \eps_k$, else set $\eps_{k+1} = \theta \eps_k$

    \item[]\makebox[0.7cm][l]{}{\bf Set} $\bar \mu_{k+1}$ by~\eqref{update_molt_2}.

\item[] {\bf End for}
\end{algorithmic}
\end{algorithm}

The next results shows that a KKT point is obtained, as a limit point, whenever we accept the Newton direction for an infinite number of iterations.

\begin{proposition}\label{prop:dk}
Let $\{x_k\}$ be a sequence generated by the Primal-Dual Augmented Lagrangian Method and let $\{x_k\}_K$ be a subsequence such that $\tilde x_k$ is accepted (i.e., $d_k$ is computed and $\|(d_k,(\tilde x_k-x_k)_{{\cal B}_k})\| \le \Delta_k$) for infinitely many iterations $k \in K$ and
\[
\lim_{k \to \infty, \, k \in K} x_{k+1} = x^*.
\]
Then, $x^*$ is a KKT point.
\end{proposition}

{\it Proof}
Since $\{\bar \mu_k\}$ is a bounded sequence and ${\cal L}_k$, ${\cal U}_k$, ${\cal N}_k$ are subsets of a finite set of indices, without loss of generality we can assume that
$\lim_{k \to \infty, \, k \in K} \bar \mu_{k+1} = \mu^*$,
${\cal L}_k = \cal L$, ${\cal U}_k = \cal U$ and ${\cal N}_k = \cal N$ (passing into a further subsequence if necessary).
Moreover, since $d_k$ is accepted for infinitely many iterations $k \in K$, without loss of generality we can also assume that $d_k$ is accepted for all $k \in K$
(passing again into a further subsequence if necessary).

Since the projection is non-expansive, for all $k \in K$ we have
\[
\|(x_{k+1},\bar \mu_{k+1}) - (x_k,\bar \mu_k)\| \le \Delta_k.
\]
Moreover, since $\Delta_{k+1} = \beta \Delta_k$, with $\beta \in (0,1)$, for all $k \in K$,
\begin{equation}\label{Delta_k_to_zero}
\lim_{k \to \infty} \Delta_k = 0.
\end{equation}
and
\[
\lim_{\substack{k \to \infty \\ k \in K}} \|x_{k+1}-x_k\| = 0.
\]
Then,
\begin{equation}\label{lim_xk}
\lim_{\substack{k \to \infty \\ k \in K}} x_k =
\lim_{\substack{k \to \infty \\ k \in K}} x_{k+1} = x^*.
\end{equation}
Since $\|d_k\| \le \Delta_k$ for all $k \in K$, from~\eqref{Delta_k_to_zero} we also have that
\begin{equation}\label{dk_to_zero}
\lim_{\substack{k \to \infty \\ k \in K}} \|d_k\| = 0.
\end{equation}

Using again the fact that the Newton direction is accepted at every iteration $k \in K$, we can write
\begin{equation}\label{newton_system}
 \left\|\begin{pmatrix}\nabla_{\cal N}L_k \\ h_k\end{pmatrix}\right\| =
 \left\|\begin{pmatrix}\nabla^2_{\cal N} L_k & \nabla_{\cal N} h_k \\ \nabla_{\cal N} h_k^T & 0\end{pmatrix}
 \begin{pmatrix}d_{{\cal N}} \\ d_{\mu}\end{pmatrix}\right\|
 \le
 \left\|\begin{pmatrix}\nabla^2_{\cal N} L_k & \nabla_{\cal N} h_k \\ \nabla_{\cal N} h_k^T & 0\end{pmatrix}\right\| \left\|
 \begin{pmatrix}d_{{\cal N}} \\ d_{\mu}\end{pmatrix}\right\|.
\end{equation}
Taking the limits for $k \to \infty$, $k \in K$,
and using~\eqref{lim_xk}, we have
\[
\begin{split}
& \lim_{\substack{k \to \infty, \\ k \in K}} \nabla_{\cal N}L_k = \lim_{\substack{k \to \infty, \\ k \in K}} \nabla_{\cal N}L_{k+1} = \nabla_{\cal N}L_(x^*,\mu^*)
\quad \text{and} \\ 
& \lim_{\substack{k \to \infty, \\ k \in K}} h_k = \lim_{\substack{k \to \infty, \\ k \in K}} h_{k+1} = h(x^*).
\end{split}
\]
Taking into account~\eqref{dk_to_zero} and~\eqref{newton_system}, we can write
\[
\nabla_{\cal N} L(x^*,\mu^*) = 0
\quad \text{and} \quad h(x^*) = 0.
\]

To conclude the proof, we have to show that
the KKT conditions are satisfied with respect to $\nabla_{x_{\cal L}} L(x^*,\mu^*)$ and $\nabla_{x_{\cal U}} L(x^*,\mu^*)$ as well.
From the instructions of the algorithm,
$(x_{k+1})_{\cal L} = (\tilde x_k)_{\cal L} = \ell_{\cal L}$ and $(x_{k+1})_{\cal U} = (\tilde x_k)_{\cal U} = u_{\cal U}$ for all $k \in K$.
Consequently,
\[
x^*_i =
\begin{cases}
\ell_i, & \quad \text{if } i \in \cal L, \\
u_i, & \quad \text{if } i \in \cal U.
\end{cases}
\]
So, using~\eqref{KKT_conditions2}, KKT conditions with respect to $\nabla_{x_{\cal L}} L(x^*,\mu^*)$ and $\nabla_{x_{\cal U}} L(x^*,\mu^*)$
hold if and only if
\begin{equation}\label{stationarity_as}
\nabla_{x_i} L(x^*,\mu^*)
\begin{cases}
\ge 0, & \quad \text{if } i \in \cal L, \\
\le 0, & \quad \text{if } i \in \cal U.
\end{cases}
\end{equation}
For any index $i\in\cal L$, from the active-set estimate~\eqref{L_as} we have
$0 \ge (d_k)_i = \ell_i - (x_k)_i \ge -\nu\sigma_i(x_k,\bar\mu_k)$
and, using the definition of $\sigma_i(x,\mu)$ given in~\eqref{def:sigma}, we get
\[
\nabla_{x_i} L_k \ge -\frac{(\ell_i-(x_k)_i)^2 + (u_i-(x_k)_i)^2}{\nu (u_i-(x_k)_i)^2} \, (d_k)_i.
\]
Similarly, for any index $i\in\cal U$ we have
$0 \le (d_k)_i = u_i - (x_k)_i \le \nu\rho_i(x_k,\bar\mu_k)$
and then
\[
\nabla_{x_i} L_k \le -\frac{(\ell_i-(x_k)_i)^2 + (u_i-(x_k)_i)^2}{\nu (l_i-(x_k)_i)^2} \, (d_k)_i.
\]
Taking the limits for $k \to \infty$, $k \in K$, and using~\eqref{lim_xk}--\eqref{dk_to_zero}, we obtain~\eqref{stationarity_as}.
\qed

In the following result, we show that any limit point of the sequence $\{x_k\}$ is either feasible for Problem~\eqref{probP} or stationary for the penalty term $\|h(x)\|^2$ of the augmented Lagrangian function, measuring the violation with respect to the equality constraints.

\begin{proposition}
Let $\{x_k\}$ be a sequence generated by the Primal-Dual Augmented Lagrangian Method and let $\{x_k\}_K$ be a subsequence such that
\[
\lim_{k \to \infty, \, k \in K} x_{k+1} = x^*.
\]
The following holds:
\begin{itemize}
    \item if $\lim_{k \to \infty} \eps_k > 0$, then $x^*$ is feasible;
    \item if $\tilde x_k$ is accepted (i.e., $d_k$ is computed and $\|(d_k,(\tilde x_k-x_k)_{{\cal B}_k})\| \le \Delta_k$) for infinitely many iterations $k \in K$, then $x^*$ is feasible (indeed, it is a KKT point);
    \item in all other cases, $x^*$ is a KKT point of the problem $\min_{\ell \le x \le u}\|h(x)\|^2$.
\end{itemize}
\end{proposition}

{\it Proof}
Let us analyze the three cases separately.
\begin{itemize}
\item If $\lim_{k \to \infty} \eps_k > 0$, from the instructions of the algorithm there exists an iteration $\hat k$ such that $\eps_{k+1} = \eps_k$ for all $k\ge \hat k$.
Therefore, $\|h(x_{k+1})\|_{\infty} \le \eta \|h(x_k)\|_{\infty}$, with $\eta \in (0,1)$, for all $k \ge \hat k$, and then $\{h(x_k)\} \to 0$, implying that $x^*$ is feasible.

\item If $\tilde x_k$ is accepted for infinitely many iterations $k \in K$, from Proposition~\ref{prop:dk} we have that $x^*$ is a KKT point, and thus it is feasible.

\item In all the other cases, we want to show that
\begin{equation}
[\nabla h(x^*) h(x^*)]_i
\begin{cases}
\ge 0, & \quad \text{if } x^*_i = \ell_i, \\
= 0, & \quad \text{if } x^*_i \in (\ell_i,u_i), \\
\le 0, & \quad \text{if } x^*_i = u_i.
\end{cases}
\end{equation}
Since $\{\bar \mu_k\}$ is a bounded sequence, without loss of generality we can assume that
$\lim_{k \to \infty, \, k \in K} \bar \mu_{k+1} = \mu^*$,
(passing into a further subsequence if necessary).
Moreover, note that there exists an iteration $\hat k \in K$ such that,
for all $k \ge \hat k$, $k \in K$,
the Newton direction $d_k$ is not accepted, that is, we compute $x_{k+1}$ such that~\eqref{opt_subprob} holds.
Since $\{\tau_k\} \to 0$, it follows that
\begin{equation}\label{lim_pg_al}
\lim_{k \to \infty} \|x_{k+1} - {\cal P}_{[\ell,u]}(x_{k+1} - \nabla_x L_a(x_{k+1},\bar \mu_k;\eps_k))\|_{\infty} = 0.
\end{equation}

Now, we distinguish three subcases.
\begin{enumerate}[label=(\roman*), leftmargin=*]
    \item $x^*_i \in (\ell_i,u_i)$. Since $\{x_{k+1}\}_K \to x^*$, there exists an iteration $\hat k \in K$ such that $(x_{k+1})_i \in (\ell_i,u_i)$ for all $k \ge \hat k$, $k \in K$.
    In view of~\eqref{lim_pg_al}, it follows that
    \[
    \lim_{k \to \infty, \, k \in K} (\nabla_x L_a(x_{k+1},\bar \mu_k;\eps_k))_i = 0
    \]
    (otherwise, if it was not true, then $\limsup_{k \to \infty, \, k \in K} |(x_{k+1} - {\cal P}_{[\ell,u]}(x_{k+1} - \nabla_x L_a(x_{k+1},\bar \mu_k;\eps_k)))_i| > 0$, leading to a contradiction with~\eqref{lim_pg_al}).
    So, there exists an iteration, that we still denote by $\hat k \in K$ without loss of generality, such that $(x_{k+1} - \nabla_x L_a(x_{k+1},\bar \mu_k;\eps_k))_i \in [\ell_i,u_i]$ for all $k \ge \hat k$, $k \in K$.
    Hence, for all $k \ge \hat k$, $k \in K$, we can write
    \[
    \begin{split}
    \tau_k & \ge \|x_{k+1} - {\cal P}_{[\ell,u]}(x_{k+1} - \nabla_x L_a(x_{k+1},\bar \mu_k;\eps_k))\|_{\infty} \\
    & \ge \bigl|(x_{k+1} - {\cal P}_{[\ell,u]}(x_{k+1} - \nabla_x L_a(x_{k+1},\bar \mu_k;\eps_k)))_i\bigr| \\
    & = \bigl|(\nabla_x L_a(x_{k+1},\bar \mu_k;\eps_k))_i\bigr| \\
    & = \biggl|\biggl(\nabla f(x_{k+1}) + \nabla h(x_{k+1}) \bar \mu_k + \frac2{\eps_k} \nabla h(x_{k+1}) h(x_{k+1})\biggr)_i\biggr|.
    \end{split}
    \]
    Multiplying the first and the last term in the above chain of inequality by $\eps_k$, we get
    \[
    \eps_k \tau_k \ge | (\eps_k \nabla f(x_{k+1}) + \eps_k \nabla h(x_{k+1}) \bar \mu_k + 2 \nabla h(x_{k+1}) h(x_{k+1}))_i|, 
    \]
    for all $k \ge \hat k$, $k \in K$.
    Taking the limits in the above inequality for $k \to \infty$, $k \in K$, the left-hand side converges to zero, since both $\{\eps_k\}$ and $\{\tau_k\}$ converge to zero, while the right-hand side converges to $|(2 \nabla h(x^*) h(x^*))_i|$,
    since $\{\eps_k\} \to 0$, $\{\nabla f(x_{k+1})\}_K \to \nabla f(x^*)$, $\{\nabla h(x_{k+1})\}_K \to \nabla h (x^*)$, $\{h(x_{k+1})\}_K \to h(x^*)$ and $\{\bar \mu_k\}_K \to \mu^*$. We thus conclude that $(\nabla h(x^*) h(x^*))_i=0$.
    \item $x^*_i = \ell_i$. Since $\{x_{k+1}\}_K \to x^*$, there exists an iteration $\hat k \in K$ such that $(x_{k+1})_i \in [\ell_i,u_i)$ for all $k \ge \hat k$, $k \in K$.
    In view of~\eqref{lim_pg_al}, it follows that
    \[
    \liminf_{k \to \infty, \, k \in K} (\nabla_x L_a(x_{k+1},\bar \mu_k;\eps_k))_i \ge 0
    \]
    (otherwise, if it was not true, then $\limsup_{k \to \infty, \, k \in K} |(x_{k+1} - {\cal P}_{[\ell,u]}(x_{k+1} - \nabla_x L_a(x_{k+1},\bar \mu_k;\eps_k)))_i| > 0$, leading to a contradiction with~\eqref{lim_pg_al}).
    So, we can write
    \[
    \liminf_{k \to \infty, \, k \in K} \biggl(\nabla f(x_{k+1}) + \nabla h(x_{k+1}) \bar \mu_k + \frac2{\eps_k} \nabla h(x_{k+1}) h(x_{k+1})\biggr)_i \ge 0.
    \]
    Multiplying the terms of the above inequality by $\eps_k$,
    and taking into account that $\{\eps_k\}~\to~0$, $\{\nabla f(x_{k+1})\}_K \to \nabla f(x^*)$, $\{\nabla h(x_{k+1})\}_K \to \nabla h (x^*)$, $\{h(x_{k+1})\}_K \to h(x^*)$ and $\{\mu_k\}_K \to \mu^*$ is bounded, we get
    \[
    \begin{split}
    & \liminf_{k \to \infty, \, k \in K} (\eps_k \nabla f(x_{k+1}) + \eps_k \nabla h(x_{k+1}) \bar \mu_k + 2\nabla h(x_{k+1}) h(x_{k+1}))_i =\\
    & \qquad\qquad = 2(\nabla h(x^*) h(x^*))_i \ge 0.
    \end{split}
    \]
    \item $x^*_i = u_i$. We obtain $(\nabla h(x^*) h(x^*))_i \le 0$ using the same arguments as in the previous case.\qed
\end{enumerate}
\end{itemize}

In order to show convergence of the algorithm to KKT points, we need to point out some properties of the approximate minimizers of the augmented Lagrangian function.
In particular, in the next lemma we show that, when we cannot use the Newton direction,
the approximate minimizers of the augmented Lagrangian function computed as in~\eqref{opt_subprob}, with $\{\tau_k\} \to 0$, satisfy the conditions stated in~\cite{andreani:2008} for the solutions of the subproblems (see Step~2 of Algorithm~3.1 in~\cite{andreani:2008}).

\begin{lemma}\label{lemma:subprob}
Let $\{x_k\}$ be a sequence generated by the Primal-Dual Augmented Lagrangian Method and let $\{x_k\}_K$ be a subsequence such that
\[
\lim_{k \to \infty, \, k \in K} x_{k+1} = x^*,
\]
with $x^*$ feasible and, for all $k \in K$, either the Newton direction $d_k$ cannot be computed (i.e., system~\eqref{KKT_systemB} does not have solutions) or $\tilde x_k$ is not accepted (i.e., $\|(d_k,(\tilde x_k-x_k)_{{\cal B}_k})\| > \Delta_k$).
Then, for all $k \in K$ there exist $\tau_{k,1} \ge 0$, $\tau_{k,2} \ge 0$, $(v_k)_i$,$(w_k)_i$,$i=1,\ldots,n$, such that
\begin{gather}
\biggl\|\nabla L_a(x_{k+1},\bar \mu_k; \eps_k) + \sum_{i=1}^n ((v_k)_i-(w_k)_i) \biggr\|_{\infty} \le \tau_{k,1}, \label{opt_subprob_1} \\
(v_k)_i \ge 0, \quad (w_k)_i \ge 0 \quad \text{and} \quad \ell_i - \tau_{k,2} \le (x_{k+1})_i \le u_i + \tau_{k,2}, \quad i = 1,\ldots,n, \label{opt_subprob_2} \\
(x_{k+1})_i > \ell_i + \tau_{k,2} \, \Rightarrow \, (w_k)_i = 0, \quad i = 1,\ldots,n, \label{opt_subprob_3} \\
(x_{k+1})_i < u_i - \tau_{k,2} \, \Rightarrow \, (v_k)_i = 0, \quad i = 1,\ldots,n, \label{opt_subprob_4} \\
\lim_{k \to \infty, \, k \in K} \tau_{k,1} = \lim_{k \to \infty, \, k \in K} \tau_{k,2} = 0. \label{opt_subprob_5}
\end{gather}
\end{lemma}

{\it Proof} First, note that the conditions on $x_k$ in~\eqref{opt_subprob_2} are satisfied for any $\tau_{k,2} \ge 0$, since we maintain feasibility with respect to the constraints $\ell \le x \le u$.
Without loss of generality, we can limit to prove that
an iteration $\hat k \in K$ exists such that~\eqref{opt_subprob_1}--\eqref{opt_subprob_4}
hold for all $k \ge \hat k$, $k \in K$, and~\eqref{opt_subprob_5} is satisfied
(for the iterations $k < \hat k$, $k \in K$, we can choose arbitrary $\tau_{k,1} \ge 0$, $\tau_{k,2} \ge 0$, $(v_k)_i$,$(w_k)_i$,$i=1,\ldots,n$, with $\tau_{k,1}$ sufficiently large, satisfying~\eqref{opt_subprob_1}--\eqref{opt_subprob_4}).

From the instructions of the algorithm, at every iteration $k \in K$ we compute $x_{k+1}$ such that~\eqref{opt_subprob} holds, with $\{\tau_k\} \to 0$.
So, we can choose $\hat k$ as the first iteration such that
\begin{equation}\label{tau_k}
\tau_k < \min_{i=1,\ldots,n}\{u_i - \ell_i\}, \quad \forall k \ge \hat k.
\end{equation}

Since the index set $\{1,\ldots,n\}$ is finite, without loss of generality we can define the subsets $I_1$, $I_2$, $I_3$ and $I_4$ (passing into a further subsequence if necessary) such that:
\begin{align*}
I_1 & = \{i \colon (x_{k+1})_i \in (\ell_i,u_i) \, \forall k \in K \text{ and } x^*_i \in (\ell_i,u_i)\}, \\
I_2 & = \{i \colon (x_{k+1})_i \in (\ell_i,u_i) \, \forall k \in K \text{ and } x^*_i \in \{\ell_i,u_i\}\}, \\
I_3 & = \{i \colon (x_{k+1})_i = \ell_i \, \forall k \in K\}, \\
I_4 & = \{i \colon (x_{k+1})_i = u_i \, \forall k \in K\}.
\end{align*}
From~\eqref{opt_subprob} and~\eqref{tau_k}, for all $k \ge \hat k$, $k \in K$, we can write
\begin{align*}
|(x_{k+1} - {\cal P}_{[\ell,u]}(x_{k+1} - \nabla_x L_a(x_{k+1},\bar \mu_k;\eps_k)))_i| \le \tau_k, \quad & i \in I_1 \cup I_2, \\
\nabla_{x_i} L_a(x_{k+1},\bar \mu_k;\eps_k) \ge -\tau_k, \quad & i \in I_3, \\
\nabla_{x_i} L_a(x_{k+1},\bar \mu_k;\eps_k) \le \tau_k, \quad & i \in I_4.
\end{align*}
For every variable $(x_k)_i$ with $i \in I_1$, we also have that $${\cal P}_{[\ell,u]}(x_{k+1} - \nabla_x L_a(x_{k+1},\bar \mu_k;\eps_k)))_i = x_{k+1} - \nabla_x L_a(x_{k+1},\bar \mu_k;\eps_k)$$ for all sufficiently large $k \in K$ (this follows from the fact that $\{(x_k)_i\}_K \to x^*_i \in (\ell_i,u_i)$ and $\tau_k \to 0$) .
So, without loss of generality we can also assume that $\hat k$ is large enough to satisfy
\begin{align*}
|\nabla_x L_a(x_{k+1},\bar \mu_k;\eps_k)_i| \le \tau_k, \quad & i \in I_1, \\
|(x_{k+1} - {\cal P}_{[\ell,u]}(x_{k+1} - \nabla_x L_a(x_{k+1},\bar \mu_k;\eps_k)))_i| \le \tau_k, \quad & i \in I_2, \\
\nabla_{x_i} L_a(x_{k+1},\bar \mu_k;\eps_k) \ge -\tau_k, \quad & i \in I_3, \\
\nabla_{x_i} L_a(x_{k+1},\bar \mu_k;\eps_k) \le \tau_k, \quad & i \in I_4.
\end{align*}

Let us rewrite the quantities within the absolute value in the second inequality as follows:
\[
(x_{k+1} - {\cal P}_{[\ell,u]}(x_{k+1} - \nabla_x L_a(x_{k+1},\bar \mu_k;\eps_k)))_i =
\nabla_{x_i} L_a(x_{k+1},\bar \mu_k;\eps_k) - (y'_k)_i+ (y''_k)_i,
\]
where $(y'_k)_i, (y''_k)_i  \ge 0$ are proper scalars.
In more detail, if $p := (x_{k+1} - \nabla_x L_a(x_{k+1},\bar \mu_k;\eps_k))_i$ is in $[\ell_i,u_i]$, then $(y'_k)_i, (y''_k)_i = 0$. On the other hand, if $p - {\cal P}_{[\ell_i,u_i]}(p)< 0$, then  $(y'_k)_i > 0$ and $(y''_k)_i = 0$; otherwise, i.e., if $p - {\cal P}_{[\ell_i,u_i]}(p)> 0$, then $(y'_k)_i = 0$ and $(y''_k)_i > 0$.
Therefore, we obtain
\begin{align*}
|\nabla_{x_i} L_a(x_{k+1},\bar \mu_k;\eps_k)| \le \tau_k, \quad & i \in I_1, \\
|\nabla_{x_i} L_a(x_{k+1},\bar \mu_k;\eps_k) - (y^1_k)_i + (y^2_k)_i| \le \tau_k, \quad & i \in I_2, \\
\nabla_{x_i} L_a(x_{k+1},\bar \mu_k;\eps_k) \ge -\tau_k, \quad & i \in I_3, \\
\nabla_{x_i} L_a(x_{k+1},\bar \mu_k;\eps_k) \le \tau_k, \quad & i \in I_4.
\end{align*}
We conclude that~\eqref{opt_subprob_1}--\eqref{opt_subprob_4} hold for all $k \ge \hat k$, $k \in K$,
with
\begin{align*}
\tau_{k,1} & = \tau_k, \\
\tau_{k,2} & =
\begin{cases}
\min_{i \in I_2 }\{\min\{(x_k)_i-\ell_i,u_i-(x_k)_i\}\}, \quad & \text{if } I_2 \ne \emptyset, \\
0, \quad & \text{if } I_2 = \emptyset,
\end{cases} \\
(v_k)_i & =
\begin{cases}
0, \quad & i \in I_1\cup I_3, \\
(y''_k)_i,\quad & i\in I_2,\\
\max\{0,-\nabla_{x_i} L_a(x_{k+1},\bar \mu_k;\eps_k)\}, \quad & i \in I_4,
\end{cases} \\
(w_k)_i & =
\begin{cases}
0, \quad & i \in I_1 \cup I_4, \\
(y'_k)_i, \quad & i \in I_2, \\
\max\{0,\nabla_{x_i} L_a(x_{k+1},\bar \mu_k;\eps_k)\}, \quad & i \in I_3,
\end{cases}
\end{align*}
and, from the above definitions, also~\eqref{opt_subprob_5} is satisfied.
\qed

Combining the above results with those stated in~\cite{andreani:2008}, we can finally show the convergence of the proposed algorithm to stationary points.
In particular, as in~\cite{andreani:2008}, we use the \textit{constant positive linear dependence} (CPLD) as constraint qualification condition.

\begin{definition}
A point $x$ is said to satisfy CPLD for Problem~\eqref{probP} if
the existence of scalars $\lambda_1,\ldots,\lambda_p$, $\pi_i \ge 0$, $i \in {\cal L}(x)$, $\varphi_j \ge 0$, $j \in {\cal U}(x)$, such that $\sum_{t=1}^p \lambda_t \nabla h_t(z) - \sum_{i \in {\cal L}(x)} \pi_i e_i + \sum_{j \in {\cal U}(x)} \varphi_j e_j = 0$ implies that,
for all $z$ in a neighborhood of $x$, the vectors $\nabla h_1(z), \ldots, \nabla h_p(z)$,$-e_i$, $i \in {\cal L}(x)$, $e_j$, $j \in {\cal U}(x)$ are linearly dependent, where
${\cal L}(x) := \{i \colon x_i = \ell_i\}$, ${\cal U}(x) := \{i \colon x_i = u_i\}$ and ${\cal N}(x) := \{1,\ldots,n\} \setminus ({\cal L}(x) \cup {\cal U}(x)).$
\end{definition}
For more details on CPLD and the relations with other constraint qualification conditions, see also~\cite{qi2000constant,andreani2005relation}.

\begin{theorem}
Let $\{x_k\}$ be a sequence generated by the Primal-Dual Augmented Lagrangian Method and let $\{x_k\}_K$ be a subsequence such that
\[
\lim_{k \to \infty, \, k \in K} x_{k+1} = x^*.
\]
The following holds:
\begin{itemize}
    \item if $\tilde x_k$ is accepted (i.e., $d_k$ is computed and $\|(d_k,(\tilde x_k-x_k)_{{\cal B}_k})\| \le \Delta_k$) for infinitely many iterations $k \in K$, then $x^*$ is a KKT point;
    \item else, if $x^*$ satisfies the CPLD constraint qualification, then $x^*$ is a KKT point.
\end{itemize}
\end{theorem}

{\it Proof}
If $\tilde x_k$ is accepted for infinitely many iterations $k \in K$, then $x^*$ is a KKT point from Proposition~\ref{prop:dk}.
Else, there exists an iteration $\hat k \in K$ such that $d_k$ is not accepted for any $k \ge \hat k$, $k \in K$,
and the algorithm reduces to a classical Augmented Lagrangian method. Then, using Lemma~\ref{lemma:subprob}, the conditions stated in~\cite{andreani:2008} for the solutions of the subproblems are satisfied and the result is obtained by the same arguments given in the proof of Theorem~4.2 in~\cite{andreani:2008}.
\qed

\section{Convergence Rate Analysis}\label{sec:rate}
In this section, we analyze the convergence rate of the proposed algorithm.
We will show that, for sufficiently large iterations, the primal-dual sequence $(x_k,\bar \mu_k)$ converges to an optimal solution $(x^*,\mu^*)$ at a quadratic rate.

In the literature, standard assumptions to prove the convergence rate of an augmented Lagrangian scheme are the \textit{linear independence constraints qualification} (LICQ), the strict complementarity and the \textit{second-order sufficient condition} (SOSC).
For Problem~\eqref{probP}, let us denote by $\sigma^*$ and $\rho^*$ the KKT multipliers at $x^*$ associated to the bound constraints $x \ge \ell$ and $x \le u$, respectively, and
\[
{\cal L}^* := \{i \colon x^*_i = \ell_i\}, \quad {\cal U}^* := \{i \colon x^*_i = u_i\},
\quad {\cal N}^* := \{1,\ldots,n\} \setminus ({\cal L}^* \cup {\cal U}^*).
\]
Then,
\begin{itemize}
    \item LICQ means that the vectors $\nabla h_1(x^*),\ldots,\nabla h_p(x^*)$, $-e_i$, $i\in {\cal L}^*$, $e_j$, $j\in {\cal U}^*$, are linearly independent;
\item SOSC means that $y^T \nabla^2_{xx} L(x^*,\mu^*) y > 0$ for all $y \in T(x^*) \setminus \{0\}$,
where
\[
\begin{split}
T(x^*) := \{y \in \R^n \colon \nabla h(x^*)^T y & = 0, \\
e_i^T y & = 0, \, i \in I_0(x^*), \\
e_i^T y & \le 0, \, i \in I_1(x^*)\},
\end{split}
\]
with
$I_0(x^*) := ({\cal L}^*\cap\{i \colon \sigma_i^* > 0\}) \cup ({\cal U}^*\cap\{i \colon \rho_i^* > 0\})$ and $I_1(x^*) := ({\cal L}^* \cup {\cal U}^*) \setminus I_0(x^*)$.
\end{itemize}
Under LICQ, strict complementarity and SOSC, if the penalty parameter $\epsilon_k \to 0$, usually it is possible to show superlinear convergence rate for augmented Lagrangian methods (see, e.g.,~\cite{bertsekas2014constrained,fernandez2012local} and the references therein).
Moreover, superlinear convergence rate is proved in~\cite{fernandez2012local}, when $\epsilon_k \to 0$, even without any constraint qualification, but requiring the starting multiplier to be in a neighborhood of a KKT multiplier satisfying SOSC.

Here, quadratic convergence rate is obtained by assuming that $\mu^*_i\in[\bar \mu_{\text{min}},\bar \mu_{\text{max}}]$ for all $i = 1,\ldots,p$,
under LICQ and the \textit{strong second-order sufficient condition} (SSOSC), where the latter means that
\[
y^T \nabla^2_{xx} L(x^*,\mu^*) y > 0, \quad \forall y \in T'(x^*) \setminus \{0\},
\]
with
\[
T'(x^*) := \{y \in \R^n \colon \nabla h(x^*)^T y = 0, \,\, e_i^T y = 0, \, i \in I_0(x^*)\}.
\]
Interestingly, our results do not need the convergence of $\{\epsilon_k\}$ to $0$.

First, we state an intermediate result ensuring that, if a sequence converges to a point where the conditions for superlinear convergence rate of the Newton direction are satisfied, then the direction is eventually accepted by the algorithm.

\begin{proposition}\label{utile:superlin}
Let $\{(z_k,\bar d_k)\}$ be a sequence of vectors such that
\[
\lim_{k\to\infty} z_k = z^* \quad \text{and} \quad
\|z_{k} + \bar d_k - z^* \| \le \alpha_k\|z_k-z^*\|,
\]
with $\{\alpha_k\} \to 0$. Then, for $k$ sufficiently large,
\[
\|\bar d_k \| \le \beta^k \Delta_0,
\]
for given $\beta \in(0,1)$ and $\Delta_0>0$.
\end{proposition}

{\it Proof}. Let $\bar k$ and $\bar \alpha$ be such that, for all $k\ge\bar k$,
\begin{equation}\label{luigi-new}
 \alpha_k<\bar \alpha<\beta<1.
\end{equation}
Therefore, we can write
\begin{align*}
\|z_{k}-z^*\| & \le \bar\alpha^{k-\bar k}\|z_{\bar k}-z^*\|, \\
\|z_{k}+d_k-z^*\| & \le \bar\alpha^{k+1-\bar k}\|z_{\bar k}-z^*\|,
\end{align*}
from which we obtain:
\[
\|\bar d_k\|\le\|z_{k}+d_k-z^*\|+ \|z_{k}-z^*\|\le \bar\alpha^{k}\,\frac{(\bar \alpha+1)}{\bar\alpha^{\bar k}}\,\|z_{\bar k}-z^*\|.
\]
By using (\ref{luigi-new}), we can set
\[
\bar\alpha=\rho\beta,\qquad\rho\in (0,1).
\]
Then, we have
\begin{equation}\label{luigi-new3}
\|\bar d_k\|\le\beta^k \rho^{k}\,\frac{(\bar \alpha+1)}{\bar\alpha^{\bar k}}\,\|z_{\bar k}-z^*\|.
\end{equation}
Since $\rho\in (0,1)$, we can conclude that, for $k$ sufficiently large, it results that
\begin{equation}\label{luigi-new4}
\rho^{k}\,\frac{(\bar \alpha+1)}{\bar\alpha^{\bar k}}\,\|z_{\bar k}-z^*\|\le \Delta_0.
\end{equation}
Now, \eqref{luigi-new3} and \eqref{luigi-new4}
conclude the proof.
\qed

Finally, we are ready to show the asymptotic quadratic rate of the primal-dual sequence
$\{(x_k,\bar \mu_k)\}$, under LICQ and SSOSC, if $\mu^*_i\in[\bar \mu_{\text{min}},\bar \mu_{\text{max}}]$ for all $i = 1,\ldots,p$.

\begin{theorem}
Let $\{x_k\}$ and $\{\bar\mu_k\}$ be the sequences generated by the Primal-Dual Augmented Lagrangian Method and assume that
\[
\lim_{k \to \infty} x_k = x^*, \quad \lim_{k\to\infty}\bar\mu_k = \mu^*,
\]
with $\mu^*_i\in[\bar \mu_{\text{min}},\bar \mu_{\text{max}}]$ for all $i = 1,\ldots,p$. Also assume that the LICQ and SSOSC hold at $(x^*,\mu^*)$. Then $\{(x_k,\bar\mu_k)\}$ converges to $(x^*,\mu^*)$ with a quadratic rate asymptotically, i.e.,
\[
\left\|\begin{matrix}x_{k+1} - x^* \\
\bar\mu_{k+1} - \mu^*\end{matrix}\right\| \le
K \left\|\begin{matrix}x_k - x^* \\
\bar\mu_k - \mu^*\end{matrix}\right\|^2
\]
for all sufficiently large $k$ and some constant $K$.
\end{theorem}

{\it Proof}~ Since LICQ and SSOSC hold at $(x^*,\mu^*)$, using~\cite[Proposition~3.1]{facchinei:1995} it follows that the following matrix is invertible for $k$ sufficiently large:
\[
\begin{pmatrix}
\nabla^2 L_k & \nabla h_k & -I_{{\cal L}_k} & I_{{\cal U}_k} \\
\nabla h_k^T & 0 & 0 & 0\\
-I_{{\cal L}_k}^T & 0 & 0 & 0\\
I_{{\cal U}_k}^T & 0 & 0 & 0
\end{pmatrix},
\]
where $I_{{\cal L}_k}$ and $I_{{\cal U}_k}$ denote the submatrices obtained from the identity matrix by discarding the columns whose indices do not belong to ${\cal L}_k$ and ${\cal U}_k$, respectively.
Consequently, for all sufficiently large $k$, the Newton direction can be computed.

Let us define $\bar d_k$ as the the Newton direction $d_k = ((d_x)_k,(d_{\mu})_k)$ augmented with the components in ${\cal B}_k$. Namely,
$\bar d_k := ((\bar d_x)_k,(\bar d_{\mu})_k)$, where
\[
(\bar d_x)_k := ((d_x)_k, (\tilde x_k - x_k)_{{\cal B}_k}) \quad \text{and} \quad (\bar d_{\mu})_k = (d_{\mu})_k
\]
(by properly reordering the entries of $x_k$). We note that
\[
\left\|\begin{matrix}x_k + (\bar d_{x})_k - x^* \\
\bar \mu_k + (\bar d_\mu)_k - \mu^*\end{matrix}\right\| =
\left\|\begin{matrix}
 (x_k + (\bar d_{x})_k - x^*)_{{\cal L}_k} \\
 (x_k + (\bar d_{x})_k - x^*)_{{\cal U}_k} \\
 (x_k + (\bar d_{x})_k - x^*)_{{\cal N}_k} \\
 \bar \mu_k + (\bar d_\mu)_k - \mu^*
\end{matrix}\right\|.
\]
By the instructions of the algorithm, when a Newton direction is used, we have
\[
(x_{k+1})_{{\cal L}_k} = (\tilde x_k)_{{\cal L}_k} = (x_k + (\bar d_{x})_k)_{{\cal L}_k} = \ell_{{\cal L}_k}
\]
and
\[
(x_{k+1})_{{\cal U}_k} = (\tilde x_k)_{{\cal U}_k} = (x_k + (\bar d_{x})_k)_{{\cal U}_k} = u_{{\cal U}_k}.
\]
So, using Proposition~\ref{prop:estim}, for all sufficiently large $k$ we have that
\begin{subequations}\label{stimaesatta}
\begin{align}
(x_{k+1} - x^*)_{{\cal L}_k} & = (\ell - x^*)_{{\cal L}_k} = 0, \\
(x_{k+1} - x^*)_{{\cal U}_k} & = (u - x^*)_{{\cal U}_k} = 0,
\end{align}
\end{subequations}
and then,
\[
\left\|\begin{matrix}x_k + (\bar d_{x})_k - x^* \\
\bar \mu_k + (\bar d_\mu)_k - \mu^*\end{matrix}\right\| =
\left\|\begin{matrix}
 (x_k + (\bar d_{x})_k - x^*)_{{\cal N}_k} \\
 \bar \mu_k + (\bar d_\mu)_k - \mu^*
\end{matrix}\right\|.
\]
For all sufficiently large $k$,
by the same arguments given in the proof of ~\cite[Proposition 4]{dipillo:2000},
there exists a constant $K$ such that
\begin{equation}\label{conv_rate_proof}
\begin{split}
\left\|\begin{matrix}x_k + (\bar d_{x})_k - x^* \\
\bar \mu_k + (\bar d_\mu)_k - \mu^*\end{matrix}\right\| & =
\left\|\begin{matrix}
 (x_k + (\bar d_{x})_k - x^*)_{{\cal N}_k} \\
 \bar \mu_k + (\bar d_\mu)_k - \mu^*
\end{matrix}\right\| \le
K \left\|\begin{matrix}
 (x_k  - x^*)_{{\cal N}_k} \\
 \bar \mu_k - \mu^*
\end{matrix}\right\|^2 \\
& \le K \left\|\begin{matrix}
 x_k  - x^* \\
 \bar \mu_k - \mu^*
\end{matrix}\right\|^2.
\end{split}
\end{equation}
The above relation implies that $\bar d_k$ satisfies the assumptions of Proposition~\ref{utile:superlin} (with $z_k = (x_k,\bar \mu_k)$ and $\alpha_k = K \|(x_k,\bar \mu_k) - (x^*,\mu^*)\|$).
Since $\|d_k\| \le \|\bar d_k\|$, by the instructions of the algorithm
the Newton direction $d_k$ is accepted for all sufficiently large $k$, so that
\[
(x_{k+1})_{{\cal N}_k} = (\tilde x_k)_{{\cal N}_k} = {\cal P}_{[\ell_{{\cal N}_k},u_{{\cal N}_k}]}((x_k)_{{\cal N}_k} + \bar d_{x_{{\cal N}_k}})
\]
and
\[
(\bar \mu_{k+1})_i = \max\{\bar\mu_{\text{min}},\min\{\bar\mu_{\text{max}},(\bar\mu_k + (\bar d_{\mu})_k)_i\}\}, \quad i=1,\ldots,p.
\]
Using~\eqref{stimaesatta}, we get
\[
\begin{split}
\|x_{k+1} - x^*\| & =
\|(x_{k+1} - x^*)_{{\cal N}_k}\| =
\|{\cal P}_{[\ell_{{\cal N}_k},u_{{\cal N}_k}]}((x_k)_{{\cal N}_k} + \bar d_{x_{{\cal N}_k}}) - (x^*)_{{\cal N}_k} \|\\
& \le \|(x_k + (\bar d_{x})_k - x^*)_{{\cal N}_k}\| \le \|x_k + (\bar d_{x})_k - x^*\|,
\end{split}
\]
where the first inequality follows from the fact that the projection operator is non-expansive
and that, for all sufficiently large $k$, from Proposition~\ref{prop:estim} we have
${\cal N}^* \subseteq {\cal N}_k$, implying that $(x^*)_{{\cal N}_k} \in (\ell_{{\cal N}_k},u_{{\cal N}_k})$.
Similarly, using again the non-expansivity of the projection operator and the assumption that
$\mu^*_i\in[\bar \mu_{\text{min}},\bar \mu_{\text{max}}]$ for all $i = 1,\ldots,p$, we have
\[
\begin{split}
\|\bar\mu_{k+1} - \mu^*\| & = \| {\cal P}_{[\bar \mu_{\text{min}}\mathbf{1},\bar \mu_{\text{max}}\mathbf{1}]}(\bar\mu_k + (\bar d_{\mu})_k) - \mu^*\| \\
& \le
\|\bar\mu_k + (\bar d_{\mu})_k - \mu^*\|,
\end{split}
\]
where $\mathbf{1}$ denotes the vector of all ones (of appropriate dimensions).
Combining these relations with~\eqref{conv_rate_proof}, for all sufficiently large $k$ we obtain
\[
\left\|\begin{matrix}x_{k+1} - x^* \\
\bar\mu_{k+1} - \mu^*\end{matrix}\right\| \le
\left\|\begin{matrix}x_k + (d_{x})_k - x^* \\
\bar\mu_k + (d_\mu)_k - \mu^*\end{matrix}\right\| \le
K \left\|\begin{matrix}x_k - x^* \\
\bar\mu_k - \mu^*\end{matrix}\right\|^2,
\]
concluding the proof. \qed

\section{Numerical Experiments}\label{sec:experiments}
This section is devoted to the description of the numerical experience with the proposed algorithm and to its comparison with other algorithms publicly available. All the numerical experiments have been carried out on an Intel Xeon CPU E5-1650 v2 @ 3.50GHz with 12 cores and 64 Gb RAM.

\medskip
{\em Problem set description}. We considered a set of $362$ general constrained problems from the CUTEst collection~\cite{gould2015cutest}, with number of variables $n \in [90, 906]$ and number of general constraints (equalities and inequalities) $m \in [1,8958]$. In particular, among the whole CUTEst problems collection, we selected all constrained problems (i.e., with at least one constraint besides bound constraints on the variables) having:
\begin{enumerate}[label=(\roman*), leftmargin=*]
    \item number of variables and constraints ``{\em user modifiable}'', or
    \item number of variables ``{\em user modifiable}'' and a fixed number of constraints, or
    \item at least 100 variables.
\end{enumerate}

\begin{figure}[htb]
 \begin{center}
  \includegraphics[width=0.8\textwidth]{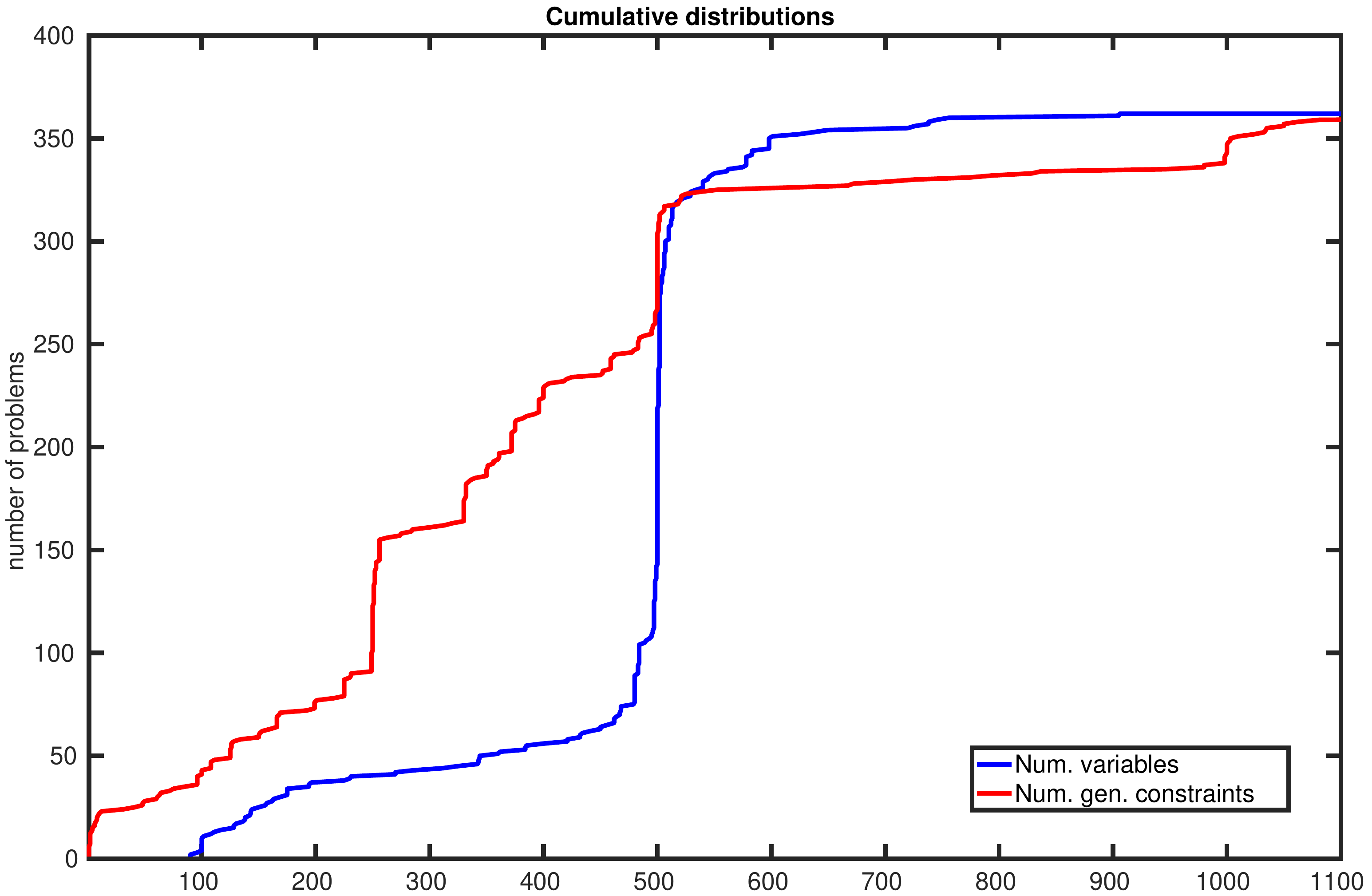}
 \end{center}
\caption{Problem set composition. The two curves represent the number of problems that have at most a given number of variables or general constraints, respectively.}\label{fig:testset}
\end{figure}

Figure~\ref{fig:testset} describes the distribution of the number of variables and number of general constraints of the considered problems.

\medskip
{\em Algorithms used in the comparison}.  We used the following algorithms:
\begin{itemize}
    \item[-] the augmented Lagrangian method implemented in the ALGENCAN (v.3.1.1) software package \cite{algencan:2008,andreani:2008};
    \item[-] the augmented Lagrangian method implemented in LANCELOT (rev.B) \cite{lancelotA,lancelotB};
    \item[-] our proposed primal-dual augmented Lagrangian method \mbox{P-D ALM} (as described in Section \ref{sec:algorithm}).
\end{itemize}

Both ALGENCAN and LANCELOT have been run using their default parameters. Note that, in its default setting, ALGENCAN uses second-order information exploiting a so-called ``acceleration strategy'', which is activated when the current primal-dual pair is sufficiently close to a KKT pair of the problem.

Our method has been implemented by modifying the code of ALGENCAN in two points:
\begin{itemize}
\item at the beginning of each iteration $k$, we inserted the computation of the active-set estimate and the Newton direction $d_k$, according to the algorithmic scheme reported in Section~\ref{sec:algorithm};
\item the approximate minimization of the augmented Lagrangian function is carried out by means of the ASA-BCP method proposed in~\cite{cristofari:2017}, in place of GENCAN~\cite{birgin2002large}.
\end{itemize}
In more detail, for every iteration $k$, in~\eqref{set_active_vars} we set $\nu = \min\{10^{-6}, \|x_k - {\cal P}_{[\ell,u]}(x_k - \nabla_x L(x_k,\bar \mu_k)) \|^{-3} \}$ and the linear system~\eqref{KKT_systemB} was solved by means of
the MA57 library~\cite{duff2004ma57}.
Note that we used the same library also in ALGENCAN.
For what concerns the inner solver ASA-BCP, it is an active-set method where, at each iteration, the variables estimated as active are set to the bounds, while those estimated as non-active are moved along a truncated-Newton direction.
In ASA-BCP, here we employed a monotone line search and, to compute the truncated-Newton direction by conjugate gradient, we used the preconditioning technique described in~\cite{birgin2014practical}, based on quasi-Newton formulas.

It is worth noticing that, in our implementation of \mbox{P-D ALM}, the test for accepting the point $\tilde x_k$ is made of two conditions, which must be both satisfied for acceptance. The first condition is that reported in Section~\ref{sec:algorithm}, i.e., $\|(d_k,(\tilde x_k-x_k)_{{\cal B}_k})\| \le \Delta_k$,
while the second condition is that $\|h(\tilde x_k)\|_{\infty} \le \eta \|h(x_k)\|_{\infty}$, i.e., the feasibility violation in $\tilde x_k$ must be sufficiently smaller than in $x_k$.
In our experience, adding this new condition leads to better results in practice.

In our experiments, for all the considered methods we used the same stopping conditions. Namely,
the algorithms were stopped when the following two conditions were both satisfied:
\begin{align*}
  \|x_k - {\cal P}_{[u,\ell]}(x_k - \nabla L(x_k,\bar\mu_k)) \|_\infty & \leq \eps_{\text{opt}} \max\{1,\|\nabla f(x_k)\|_\infty\}, \\
  \|h(x_k)\|_\infty & \leq \eps_{\text{feas}}\|h(x_0)\|_\infty,
\end{align*}
where $x_0$ is the initial point and $(x_k,\bar\mu_k)$ is the primal-dual pair at iteration $k$, with $\eps_{\text{opt}} = \eps_{\text{feas}} = 10^{-6}$.
Moreover, we inserted a maximum number of (outer) iterations equal to $400$ and a time limit of $3600$ seconds.

\begin{figure}[htb]
 \centering
  \includegraphics[width=\textwidth]{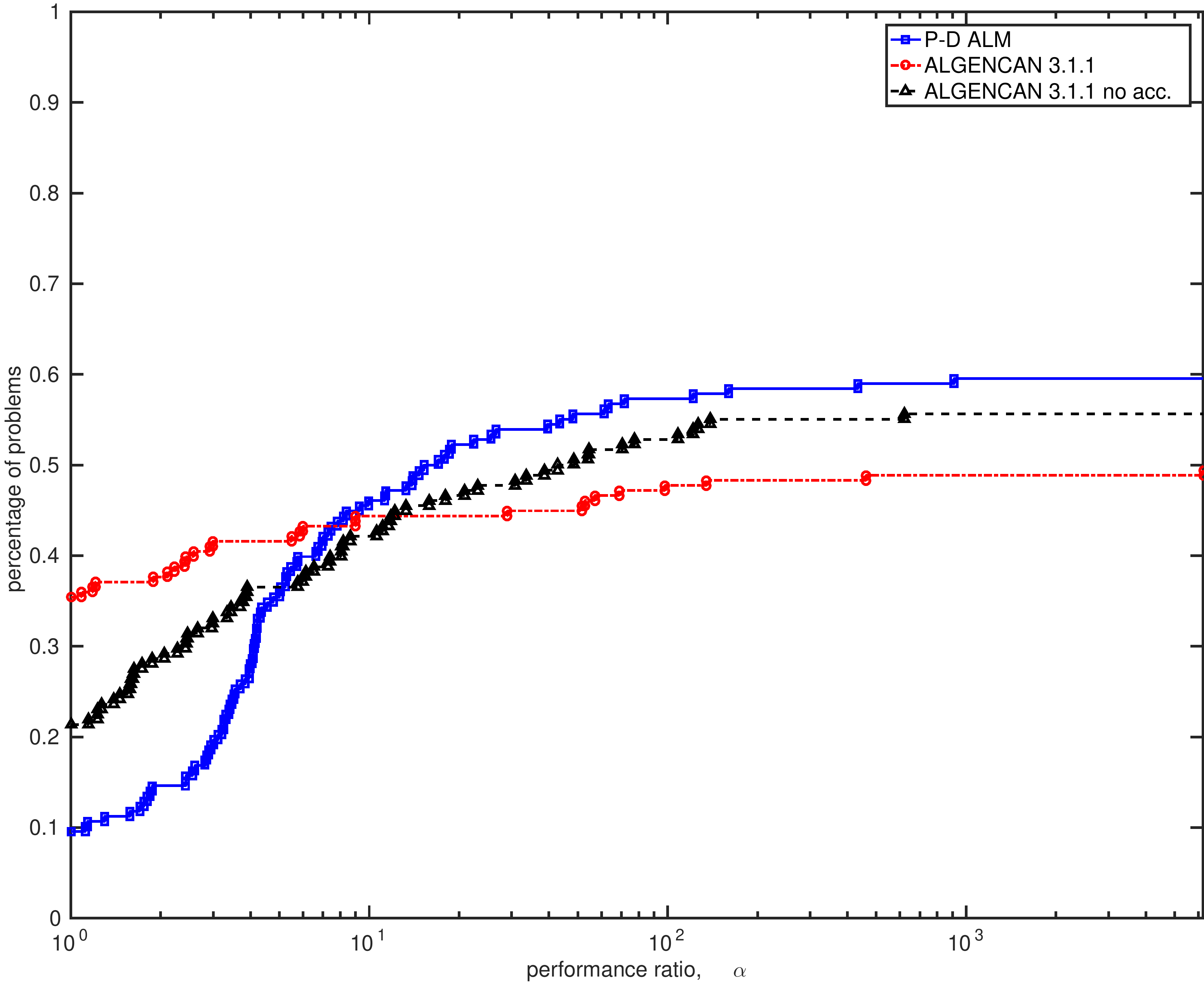}
 \caption{Comparison between \mbox{P-D ALM} and ALGENCAN with and without acceleration step, using performance  profiles with respect to CPU time. Note that ``ALGENCAN 3.1.1 no acc.'' refers to the version of ALGENCAN not using second-order information, i.e., skipping the so-called acceleration phase.}\label{fig:compatre}
\end{figure}

In Figure~\ref{fig:compatre}, we start by comparing \mbox{P-D ALM} against ALGENCAN with and without acceleration phase (note that the acceleration phase in ALGENCAN is where second-order information come into play) using the performance profiles~\cite{dolan:2002} with respect to CPU time. Note that the performance profiles are obtained on the subset of problems where at least one solver requires more than 10 seconds of CPU time. As it can be seen, ALGENCAN (using second-order information) is the most efficient solver but the least robust one. On the other hand, \mbox{P-D ALM} is considerably more robust than both the versions of ALGENCAN. One possible reason for \mbox{P-D ALM} being less efficient than ALGENCAN can be the following: in \mbox{P-D ALM} we try to use the second-order direction as much as possible, whereas second-order information is used in ALGENCAN only when the current primal-dual point is sufficiently close to a KKT pair. This could explain our larger computational times and the behaviour of the reported performance profiles.

\begin{figure}[htb]
 \centering
 \subfloat[]
    {\includegraphics[scale=0.175]{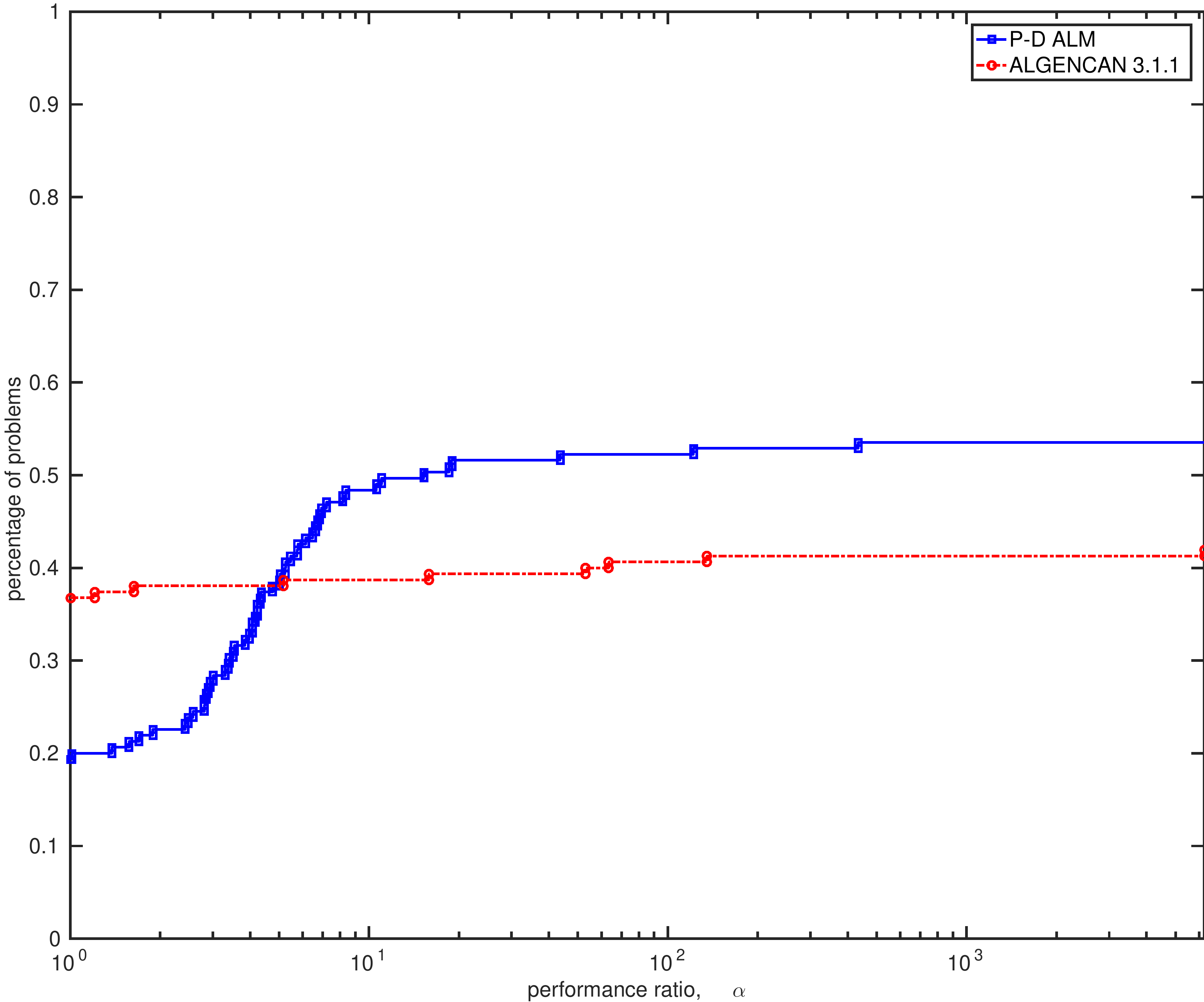}\label{fig:algencan_vs_pdalm}}
    \quad
 \subfloat[]
    {\includegraphics[scale=0.175]{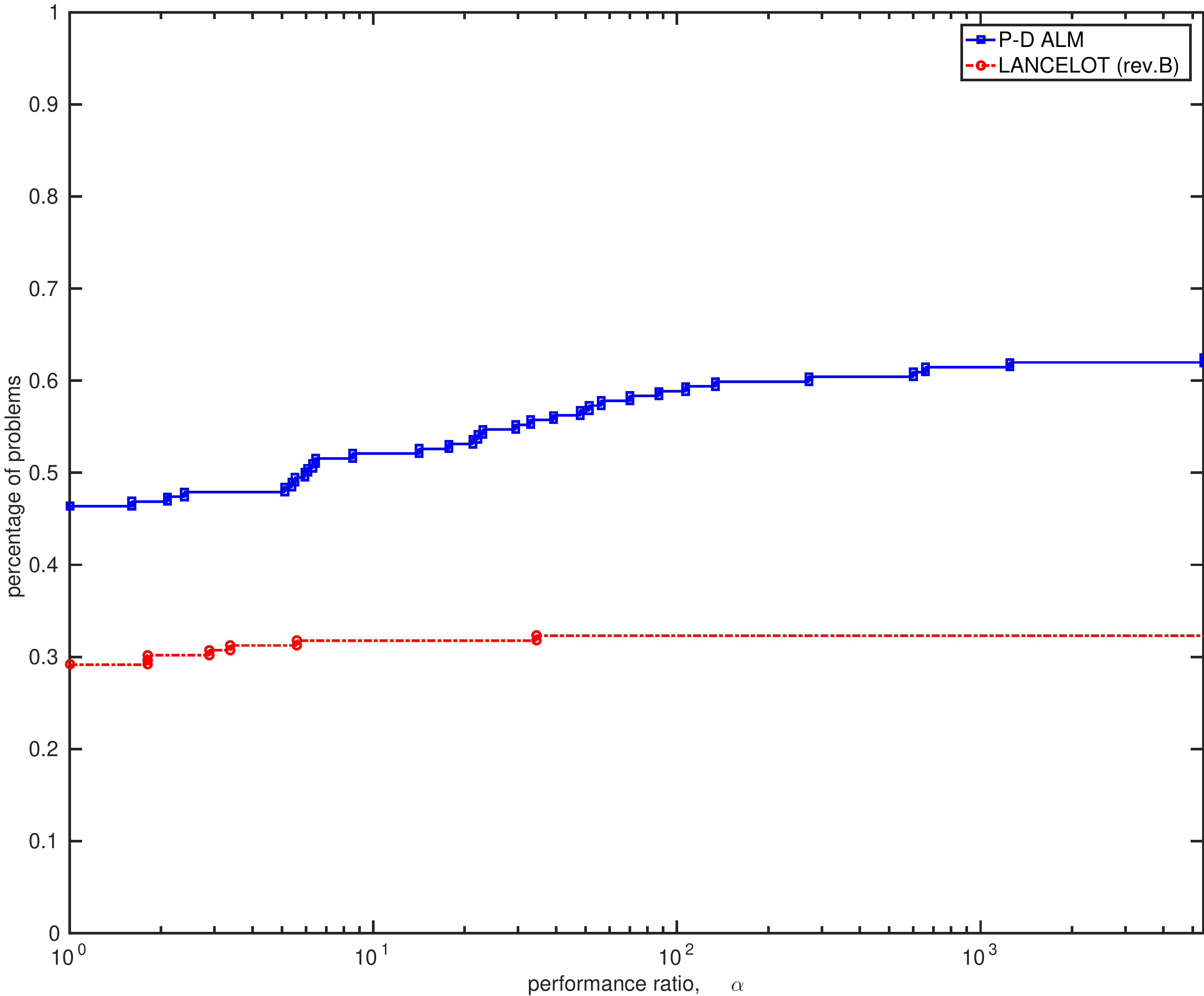}\label{fig:pdalm_vs_lancelot}}
\caption{(a) Comparison between \mbox{P-D ALM} and ALGENCAN, using performance  profiles with respect to CPU time.
 (b) Comparison between \mbox{P-D ALM} and LANCELOT, using performance  profiles with respect to CPU time.}
\end{figure}

In Figure \ref{fig:algencan_vs_pdalm}, we report the comparison between ALGENCAN and \mbox{P-D ALM}. We note that, even though ALGENCAN is slightly better than \mbox{P-D ALM} in terms of efficiency, it is outperformed by our proposed method in terms of robustness. Furthermore, we note that the two performance profiles intersect at, approximately, $\alpha \simeq 5$, i.e., both algorithms solve the same percentage of problems in at most 5 times the CPU time of the best performing solver.

In Figure \ref{fig:pdalm_vs_lancelot}, we report the comparison between \mbox{P-D ALM} and LANCELOT (rev. B). In this case, \mbox{P-D ALM} is clearly the best performing solver both in terms of efficiency and robustness.

\begin{figure}[htb]
 \centering
  \includegraphics[width=\textwidth]{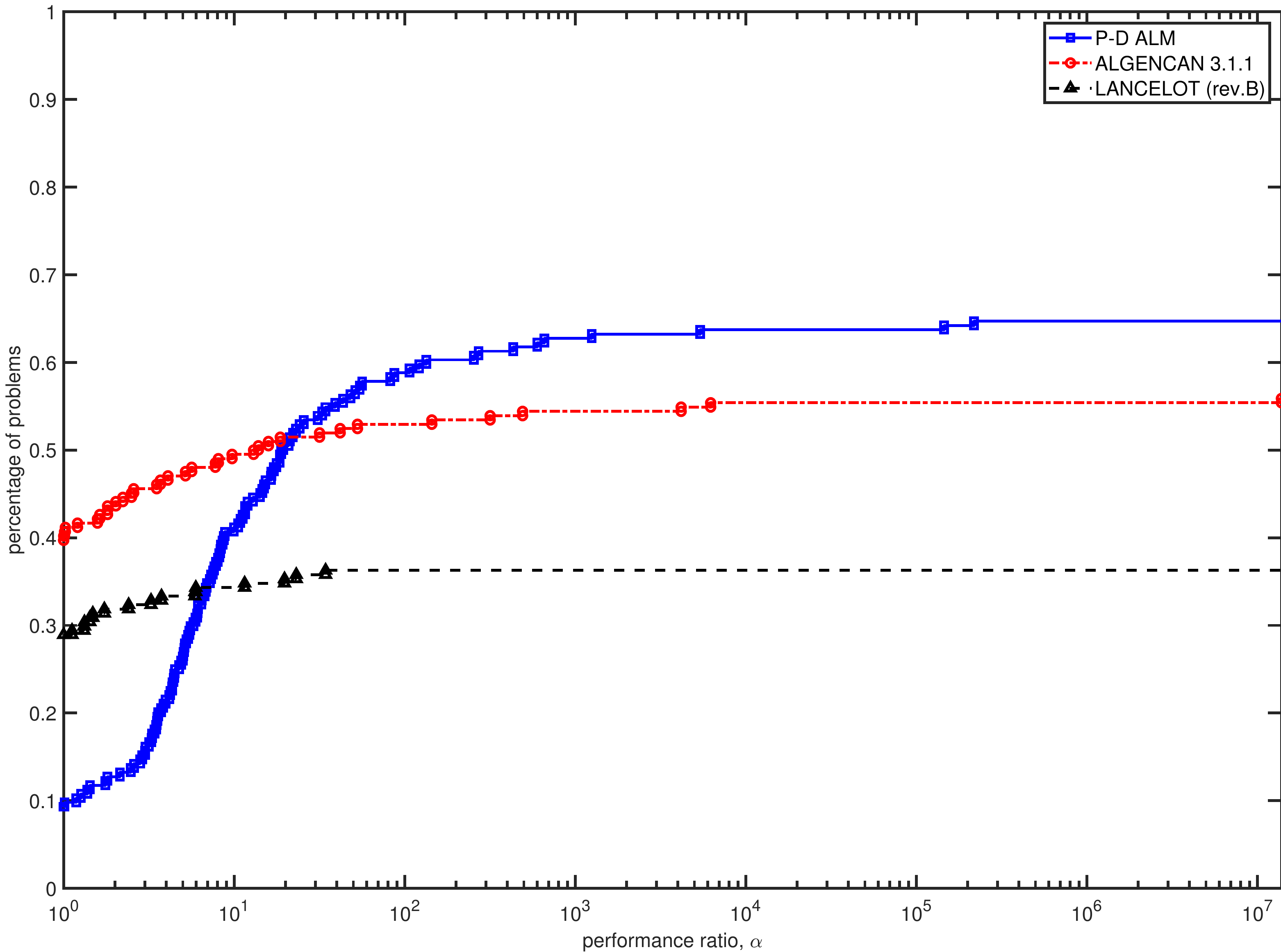}
 \caption{Comparison between \mbox{P-D ALM}, ALGENCAN and LANCELOT.}\label{fig:comp_pdalm_algencan_lancelot}
\end{figure}

Finally, we notice that ALGENCAN, LANCELOT and \mbox{P-D ALM} solve, respectively, $272$, $232$ and $290$ problems out of $362$. The comparison among the three solvers is reported in Figure~\ref{fig:comp_pdalm_algencan_lancelot}.


\section{Conclusions}\label{sec:conclusions}
In this paper, we presented a new method for nonlinear optimization problems with equality constraints and bound constraints. Starting from the augmented Lagrangian scheme implemented in ALGENCAN, we used a tailored active-set strategy to compute a Newton-type direction with respect to the variables estimated as non-active, while the variables estimated as active are set to the bounds.
If this direction satisfies a proper test, an augmented Lagrangian function is minimized by means of an efficient solver recently proposed in the literature.
We proved convergence to stationary points and, under standard assumptions, an asymptotic quadratic convergence rate.
The numerical results show the effectiveness of the proposed method.

\clearpage
\bibliography{pd_alm_2021_biblio}

\end{document}